\documentclass[english]{amsart}
\usepackage{amsmath}
\usepackage{amssymb}
\usepackage{amsthm}
\usepackage{babel}
\newtheorem{lem}{Lemma}[section]
\newtheorem{prop}{Proposition}[section]

\newtheorem{coro}{Corollary}[section]
\theoremstyle{definition}
\newtheorem{rem}{Remark}
\author{Antonio M. Oller Marc\'{e}n}
\title{The $SL(2,\mathbb{C})$ character variety of a class of torus knots}
\date{}

\begin{document}
\maketitle

\begin{abstract}
In this paper we present some families of polynomials and use them
to find, using the techniques in \cite{gma}, a defining polynomial
for the $SL(2,\mathbb{C})$ character variety (as defined in
\cite{cus}) of the torus knots of type $(m,2)$ with $m>1$ being an
odd integer.
\end{abstract}
\section{The character variety of a finitely presented group}

Let us consider a finitely presented group
$$G=\langle x_1,\dots,x_n\ |\ r_1,\dots,r_s\rangle$$
and let be $\rho:G\longrightarrow SL(2,\mathbb{C})$ be a
representation, i.e, a group homomorphism. It is clear that $\rho$
is completely determined by the $n$-tuple
$(\rho(x_1),\dots,\rho(x_n))$ and thus we can define
$$R(G)=\{(\rho(x_1),\dots,\rho(x_n))\ |\ \rho\ \textrm{is a
representation of $G$}\}\subseteq\mathbb{C}^{4n}$$ which is (see
\cite{cus}) an (up to canonical isomorphism) well-defined affine
algebraic set.

Recall that we define the character $\chi_{\rho}:G\longrightarrow F$
(see \cite{isa}) of a representation $\rho:G\longrightarrow GL(n,F)$
by $\chi_{\rho}(g)=\textrm{tr}(\rho(g))$, two representations $\rho$
and $\rho'$ having the same character if and only if they are
equivalent; i.e, if there exists $P\in GL(n,F)$ such that
$\rho'(g)=P^{-1}\rho(g)P$ for all $g\in G$. Now choose any $g\in G$
and define $\tau_g:R(G)\longrightarrow\mathbb{C}$ by
$\tau_g(\rho)=\chi_{\rho}(g)$. It is easily seen that $T=\{\tau_g\
|\ g\in G\}$ is a finitely generated ring (\cite{cus} Proposition
1.4.1.) and moreover it can be shown using some identities holding
in $SL(2,\mathbb{C})$ (see \cite{gma} Corollary 4.1.2.) that $T$ is
generated by the set:
$$\{\tau_{x_i},\tau_{x_jx_k},\tau_{x_lx_mx_p}\ |\ 1\leq i\leq n,\
1\leq j<k\leq n,\ 1\leq l<m<p\leq n\}$$

Now choose $\gamma_1,\dots\gamma_{\nu}\in G$ such that $T=\langle
\tau_{\gamma_i}\ |\ 1\leq i\leq\nu\rangle$ and define the map
$t:R(G)\longrightarrow\mathbb{C}^{\nu}$ by
$t(\rho)=(\tau_{\gamma_1}(\rho),\dots,\tau_{\gamma_\nu}(\rho))$. Put
$X(G)=t(R(G))$, then $X(G)$ is an algebraic variety (\cite{cus}
Corollary 1.4.5.) which is well-defined up to canonical isomorphism
and is called the $SL(2,\mathbb{C})$ character variety of the group
$G$. Observe that $\displaystyle{\nu=\frac{n(n^2+5)}{6}}$.

For every $0\leq j\leq n$ and for every $1\leq i\leq s$ we have that
$\tau_{r_ix_j}-\tau_{x_j}=p_{ij}$ is a polynomial with rational
coefficients in the variables $\{\tau_{x_{i_1}\dots x_{i_m}}\ |\
m\leq3\}$. With this definition we have that (see \cite{gma} Theorem
3.2.)
$$X(G)=\{\overline{x}\in\mathbb{C}^{\nu}\ |\
p_{ij}(\overline{x})=0,\ \forall i,j\}$$

\section{Torus knots}
Recall that $\mathbb{R}^2$ is the universal covering of the torus
$T^2$. We define the action
$\phi:(\mathbb{Z}\times\mathbb{Z})\times\mathbb{R}^2\longrightarrow\mathbb{R}^2$
by $\phi((m,n),(x,y))=(x+m,y+n)$ and we have that
$\mathbb{R}^2/(\mathbb{Z}\times\mathbb{Z})\substack{\cong\\
\mu}T^2$. If we now take the family $\{r_p:y=px\ |\
p\in\mathbb{R}\}$ of straight lines passing through the origin, it
is easily seen that if $p$ is irrational then $\mu(r_p)$ is dense in
$T^2$ and if $\displaystyle{p=\frac{m}{n}}$ then $\mu(r_p)\subseteq
T^2\subseteq\mathbb{R}^3$ is a knot. We denote this knot by
$K_{\frac{m}{n}}$ and call it the torus knot of type (m,n) (see
\cite{rol} Chapter 3 for further considerations).

If we denote, as usual, by $G(K)$ the fundamental group of the
exterior of any knot $K$ we can see that
$$G(K_{\frac{m}{n}})\cong\langle A,B\ |\ A^m=B^n\rangle$$
Now let us define the following group:
$$H_m=\langle x,y\ |\ \underbrace{xyxy\dots yx}_{\textrm{length $m$}}=\overbrace{yxyx\dots
xy}^{\textrm{length $m$}}\rangle$$ Then we have
\begin{lem}\label{toro}
Let $m\geq1$ be an odd integer. Then $G(K_{\frac{m}{2}})\cong H_m$.
\end{lem}
\begin{proof}
We define $\varphi:H_m\longrightarrow G(K_{\frac{m}{2}})$ by
$\varphi(x)=B^{-1}A^{\frac{m+1}{2}}$ and
$\varphi(y)=A^{-\frac{m-1}{2}}B$. On the other hand, define
$\psi:G(K_{\frac{m}{2}})\longrightarrow H_m$ by $\psi(A)=yx$ and
$\psi(B)=\overbrace{yxyx\dots y}^{\textrm{length $m$}}$. Seeing that
these homomorphisms are well defined and are each other's inverse is
straightforward.
\end{proof}

\section{Some families of polynomials}
We will start this section by defining recursively the following
family of polynomials:
$$q_1(T)=T-2$$
$$q_2(T)=T+2$$
$$\prod_{1\neq d|n} q_d\left(X+\frac{1}{X}\right)=\frac{X^{n-1}+X^{n-2}+\dots
+X+1}{X^{\frac{n-1}{2}}}\quad \textrm{if $n$ is odd}$$
$$\prod_{1,2\neq d|n} q_d\left(X+\frac{1}{X}\right)=\frac{X^{n-2}+X^{n-4}+\dots
+X^2+1}{X^{\frac{n-2}{2}}}\quad \textrm{if $n$ is even}$$
\begin{rem}
If we recall the recursive definition of the cyclotomic polynomials
(see \cite{hun} Chapter 5) by
$$\prod_{d|n}g_d(T)=T^n-1$$
then it is easily seen that
$\displaystyle{g_r(X)=X^{\frac{\varphi(r)}{2}}q_r\left(X+\frac{1}{X}\right)}$
where $\varphi$ is the Euler function.
\end{rem}
Now we introduce another family of polynomials:
$$p_1(X)=X$$
$$p_2(X)=X^2-2$$
$$p_n(X)=Xp_{n-1}(X)-p_{n-2}(X),\ \forall n\geq3$$
\begin{rem}
Let $G$ be a group and $\rho:G\longrightarrow SL(2,\mathbb{C})$ a
representation. Then $p_n(\textrm{tr}\rho(x))=\textrm{tr}\rho(x^n)$
for every $n\geq1$. For the sake of completeness we will set, where
necessary, $p_0(X)=1$.
\end{rem}
We have the following relationship between the families we have just
defined:
\begin{prop}
$$p_n(X)-2=q_1(X)\prod_{1\neq d|n}q_d^2(X)\quad \textrm{if $n$ is
odd}$$
$$p_n(X)-2=q_1(X)q_2(X)\prod_{1,2\neq d|n}q_d^2(X)\quad \textrm{if $n$ is even}$$
\end{prop}
\begin{proof}
We will just show the result for an odd $n$, the even case being
completely analogous.

Consider the cyclic group $G=\langle x\rangle$ and a representation
$\rho:G\longrightarrow SL(2,\mathbb{C})$. We can suppose,
conjugating if necessary, that $\rho(x)=\begin{pmatrix}a & b \\ 0 &
a^{-1}\end{pmatrix}$. In such a case it must be
$\rho(x^n)=\rho(x)^n=\begin{pmatrix}a^n & c \\ 0 &
a^{-n}\end{pmatrix}$

Set $X=\textrm{tr}(\rho(x))=a+a^{-1}$, then
\begin{equation*}\begin{split}p_n(X)-2&=\textrm{tr}(\rho(x^n))-2=a^n+a^{-n}-2=
\frac{(a^n-1)^2}{a^n}=\frac{1}{a^n} \left(\prod_{d|n}g_d(a)\right)^2
\\&=\frac{(a-1)^2}{a^n}\left(\prod_{1\neq
d|n}g_d(a)\right)^2\\&=\frac{(a+a^{-1}-2)a}{a^n}\left(\prod_{1\neq
d|n}a^{\frac{\varphi(d)}{2}}q_d(a+a^{-1})\right)^2\\&=q_1(X)\prod_{1\neq
d|n}q_d^2(X)
\end{split}\end{equation*}
where the identity $\displaystyle{\sum_{d|n}\varphi(d)=n}$ was used.
\end{proof}
\begin{rem}
The roots of $p_n(X)-2$ are precisely the possible values of
$\textrm{tr}(\rho(x))$ if $\rho:G\longrightarrow SL(2,\mathbb{C})$
is a representation and $x^n=1$.
\end{rem}

Let $R$ be any ring and take
$g(T)=\displaystyle{\sum_{i=0}^na_iT^i}\in R[T]$. We define
$*:R[T]\longrightarrow R[T]$ by
$g^{*}(T)=\displaystyle{\sum_{i=0}^n(-1)^{n-i}a_iT^i}$. In the next
lemma we show some interesting properties of this application.

\begin{lem}\label{ast}
Given $g,h\in R[T]$ we have:
\begin{itemize}
\item[a)] $g^{**}=g$.
\item[b)] $(gh)^{*}=g^{*}h^{*}$.
\item[c)] If $g(T)=\displaystyle{\sum_{i=0}^na_iT^i}$, then $g^{*}=g$ if
and only if $a_i=0$ for every $i$ such that $(n-i)\equiv 1\ (mod\
2)$.
\end{itemize}
\end{lem}
\begin{proof}
c) is trivial. a) and b) follow from the identity
$g^{*}(T)=(-1)^{\textrm{deg}g}g(-T)$.
\end{proof}
We can use the involution just defined to show another relation
between our two families of polynomials.
\begin{prop}\label{prop}
If $s\geq 1$ is an integer, then
$$\sum_{i=0}^s(-1)^{i}p_{s-i}(Z)=\prod_{1\neq
d|2s+1}q_d^{\ast}(Z)$$
\end{prop}
\begin{proof}
We observe that the degree of every term in $p_s(Z)$ has the same
parity as $s=\textrm{deg}\ p_s(Z)$. This fact together with the
definition of $*$ shows that
$$\left(\sum_{i=0}^s(-1)^{i}p_{s-i}(Z)\right)^{\ast}=\sum_{i=0}^{s}p_i(Z)$$
Now, we claim that
$$\sum_{i=0}^{s}p_i(Z)=\prod_{1\neq d|2s+1}q_d(Z)$$

We will proof this by induction on $s$, the case $s=1$ being trivial
as $p_0(Z)+p_1(Z)=1+Z=q_3(Z)$. Now let $s>1$ be an odd integer (the
even case is similar), by hypothesis we have
$$\sum_{i=0}^sp_i(Z)=\sum_{i=0}^{s-1}p_i(Z)+p_s(Z)=\prod_{1\neq
d|2s-1}q_d(Z)+p_s(Z)$$ and thus, setting
$\displaystyle{Z=X+\frac{1}{X}}$ one obtains:
\begin{equation*}\begin{split}\sum_{i=0}^sp_i\left(X+\frac{1}{X}\right)&=\prod_{1\neq
d|2s-1}q_d\left(X+\frac{1}{X}\right)+p_s\left(X+\frac{1}{X}\right)\\&=\frac{\displaystyle{\sum_{i=0}^{2s-2}X^{i}}}{X^{s-1}}+q_1\left(X+\frac{1}{X}\right)\prod_{1\neq
d|s}q_s^2\left(X+\frac{1}{X}\right)+2\\&=\frac{\displaystyle{\sum_{i=0}^{2s-2}X^{i}}}{X^{s-1}}+\frac{(X-1)^2}{X}\frac{\left(\displaystyle{\sum_{i=0}^{s-1}X^{i}}\right)^2}{X^{s-1}}+2\\&=
\frac{\displaystyle{\sum_{i=0}^{2s-2}X^{i}}}{X^{s-1}}+\frac{X^{2s}+1}{X^s}=\frac{\displaystyle{\sum_{i=0}^{2s}X^i}}{X^s}=
\prod_{1\neq
d|2s+1}q_d\left(X+\frac{1}{X}\right)\end{split}\end{equation*}

The proof is now completed by applying \ref{ast} a), b).
\end{proof}

\section{The $SL(2,\mathbb{C})$ character variety of the knots
$K_{\frac{m}{2}}$}

The objective of this section is to give a generating family of
polynomials for $X(G)$ with $G=G(K_{\frac{m}{2}})$ with $m>1$ an odd
integer. In \ref{toro} we shew the isomorphism
$G(K_{\frac{m}{2}})\cong H_m$ so it is enough to find such a family
for $X(H_m)$.

Before going into our main result we have to introduce another
polynomial. We set $h(X,Z)=X^2-Z$ and $k(X)=X^2-2$. Now we define
$$\alpha_{l}(X,Z)=\begin{cases} h(X,Z) & \textrm{if $l$ is even}\\
k(X) & \textrm{if $l$ is odd}\end{cases}$$ and finally we write for
$s\geq 1$
$$
f_s(X,Z)=p_s(Z)(h(X,Z)-1)+\sum_{i=1}^{s}(-1)^{i}p_{s-i}(Z)\alpha_{i}(X,Z)$$
With these definitions we can prove the following
\begin{prop}\label{prop1}
If $m>1$ is and odd integer, then
$$
X(H_m)=\{(X,Z)\in\mathbb{C}^2\ |\ f_{\frac{m-1}{2}}(X,Z)=0\}$$
\end{prop}
\begin{proof}
We set $w=\underbrace{xyxy\dots yx}_{length\
m}\underbrace{y^{-1}x^{-1}y^{-1}x^{-1}\dots y^{-1}}_{length\ m}$.
Then, using Theorem 3.2 in \cite{gma}, we have
$$X(H_m)=\{(X,Y,Z)\in\mathbb{C}^3\ |\
p_0(X,Y,Z)=p_1(X,Y,Z)=p_2(X,Y,Z)=0\}$$ where
$$
X=\tau_x  \qquad  p_0(X,Y,Z)=\tau_w-\tau_1 $$ $$Y=\tau_y  \qquad
p_1(X,Y,Z)=\tau_{wx}-\tau_x $$ $$Z=\tau_{xy}  \qquad
p_2(X,Y,Z)=\tau_{wy}-\tau_y
$$

Now, $wy=\underbrace{xyx\dots y}_{length\ m-1}x(\underbrace{xyx\dots
y}_{length\ m-1})^{-1}$ so we have $\tau_{wy}=\tau_x$ obtaining that
$p_2(X,Y,Z)=X-Y$.

On the other hand $\tau_{wx}=\tau_w\tau_x-\tau_{wx^{-1}}$ and
$wx^{-1}=\underbrace{xy\dots x}_{length\
m}y^{-1}(\underbrace{xy\dots x}_{length\ m})^{-1}$ so we get
$\tau_{wx^{-1}}=\tau_{y^{-1}}=\tau_y$ and thus
$p_1(X,Y,Z)=\tau_{wx}-\tau_x=\tau_w\tau_x-\tau_y-\tau_x=\tau_x(\tau_w-1)-\tau_y=Xp_0(X,Y,Z)+X-Y$.

Set now $w_1=(xy)^{\frac{m-1}{2}}$ and
$w_2=(yx)^{\frac{m-1}{2}}yx^{-1}$. Then it is easy to see that
$p_0(X,Y,Z)=\tau_w$ vanishes if and only if
$f(X,Y,Z)=\tau_{w_2}-\tau_{w_1}$ does. Let us compute now this
polynomial.

Firstly it is obvious by definition that
$\tau_{w_1}=p_{\frac{m-1}{2}}(Z)$. In addition we have
$\tau_{w_2}=\tau_{(yx)^{\frac{m-1}{2}}}\tau_{yx^{-1}}-
\tau_{(xy)^{\frac{m-3}{2}}xx}=p_{\frac{m-1}{2}}(Z)(XY-Z)-\tau_{(xy)^{\frac{m-3}{2}}xx}$.
Moreover we see that
$\tau_{(xy)^{\frac{m-3}{2}}xx}=\tau_{(xy)^{\frac{m-3}{2}}}\tau_{x^2}-\tau_{(yx)^{\frac{m-5}{2}}yx^{-2}}=p_{\frac{m-3}{2}}(Z)(X^2-2)-\tau_{(yx)^{\frac{m-5}{2}}yx^{-2}}$
so it is enough to iterate the process.

By now we have obtained
\begin{equation*}\begin{split}X(H_m)&=\{(X,Y,Z)\in\mathbb{C}^3\ |\
f(X,Y,Z)=0=X-Y\}\\ &\cong\{(X,Z)\in\mathbb{C}^2\ |\
f(X,X,Z)=0\}\end{split}\end{equation*} and this completes the proof
as the equality $f(X,X,Z)=f_{\frac{m-1}{2}}(X,Z)$ is just a
straightforward computation.
\end{proof}

Let us rewrite now the polynomial $f_s(X,Z)$ in a different way. In
fact we can see that
\begin{equation*}\begin{split}f_s(X,Z)=&(X^2-Z-2)\left(\sum_{i=0}^s(-1)^{i}p_{s-i}(Z)\right)\\&+p_s(Z)+\sum_{i=1}^s(-1)^{i}\beta_{i}(Z)p_{s-i}(Z)\end{split}\end{equation*}
where $\beta_k(Z)=\begin{cases} Z & \textrm{if $k$ is odd}\\
2 & \textrm{if $k$ is even}\end{cases}$

\begin{lem}\label{lema}
$p_s(Z)+\displaystyle{\sum_{i=1}^s(-1)^{i}\beta_{i}(Z)p_{s-i}(Z)}=0$
\end{lem}
\begin{proof}
It is enough to use the fact that $p_s(Z)-Zp_{s-1}(Z)=-p_{s-2}(Z)$.
\end{proof}

\begin{coro}
If $m>1$ is an odd integer, then
$$X(H_m)\cong\{(X,Z)\in\mathbb{C}^2\ |\ (X^2-Z-2)\prod_{1\neq
d|m}q_d^{\ast}(Z)=0\}$$
\end{coro}
\begin{proof}
Just apply Proposition \ref{prop} and  Lemma \ref{lema} to
Proposition \ref{prop1}
\end{proof}
\begin{lem}
Let
$\{a_1,\overline{a}_1,\dots,a_{\frac{\varphi(r)}{2}},\overline{a}_{\frac{\varphi(r)}{2}}\}$
be set of the $\varphi(r)$ primitive $r$th roots of unity. Then
$$q_r(Z)=\prod_{i=1}^{\frac{\varphi(r)}{2}}(Z-2\textrm{Re}(a_i))$$
\end{lem}
\begin{proof}
Recall that, for $r>2$ we have
$\displaystyle{g_r(X)=X^{\frac{\varphi(r)}{2}}q_r\left(
X+\frac{1}{X}\right)}$ with $g_r(X)$ being the $r$th cyclotomic
polynomial. As for all $1\leq j\leq \frac{\varphi(r)}{2}$ it holds
that $\frac{1}{a_j}=\overline{a_j}$ we obtain that $q_r(Z)$ has
exactly $\frac{\varphi(r)}{2}$ different roots, namely
$\{2\textrm{Re}(a_1),\dots,2\textrm{Re}(a_{\frac{\varphi(r)}{2}})\}$.
This together with the fact that the degree of $q_r(Z)$ is
$\frac{\varphi(r)}{2}$ completes the proof.
\end{proof}

This lemma allows us to go one step further in our description of
the curve $X(H_m)$.

\begin{coro}
Let $m>1$ be an odd integer. In the complex plane $(X,Z)$ the curve
$X(H_m)$ consists of the parabola $Z=X^2-2$ and the union of
$\frac{m-1}{2}$ horizontal lines of the form $Z=-2\textrm{Re}(w)$,
being $1\neq w$ an $m$th root of unity.
\end{coro}
\begin{proof}
It is enough to apply the previous lemma together with the fact that
given a polynomial $g$, then a number $a$ is a root of $g$ if and
only if $-a$ is a root of $g^{*}$.
\end{proof}

\end{document}